\documentclass{amsart}
\usepackage{amssymb}
\usepackage[abbrev]{amsrefs}
\numberwithin{equation}{section}
\def\N{\mathbb N}
\providecommand{\norm}[1]{\left\lVert#1\right\rVert}
\providecommand{\abs}[1]{\left\lvert#1\right\rvert}
\newcommand{\ip}[2]{\left\langle #1, #2 \right\rangle}
\newcommand{\VI}{\operatorname{VI}}
\newcommand{\Fix}{\operatorname{Fix}}
\theoremstyle{plain}
\newtheorem{theorem}{Theorem}[section]
\newtheorem{lemma}[theorem]{Lemma}
\newtheorem{problem}[theorem]{Problem}
\theoremstyle{remark}
\newtheorem{remark}[theorem]{Remark}

\allowdisplaybreaks

\title[A note on the hybrid steepest descent methods]
{A note on the hybrid steepest descent methods}
\date{\today}

\author{Koji~Aoyama}
\address[Koji~Aoyama]{Department of Economics, Chiba University, 
Yayoi-cho, Inage-ku, Chi\-ba-shi, Chiba 263-8522, Japan}
\email{aoyama@le.chiba-u.ac.jp}

\author{Yasunori~Kimura}
\address[Yasunori~Kimura]{Department of Mathematical and Computing Sciences,
Tokyo Institute of Technology,
Ookayama, Meguro-ku, Tokyo 152-8552, Japan}
\email{yasunori@is.titech.ac.jp}



\begin{document}

\begin{abstract}
 The aim of this paper is to prove that, in an appropriate setting,  
 every iterative sequence generated by the hybrid steepest descent
 method is convergent whenever so is every iterative sequence generated
 by the Halpern type iterative method. 
\end{abstract}

\maketitle

\section{Introduction}
In this paper, we consider the following variational inequality problem
in a real Hilbert space $H$: Find $z\in F$ such that 
\[
 \ip{y-z}{Az}\geq 0 \text{ for all }y\in F,
\]
where $F$ is the set of common fixed points of a sequence $\{T_n \}$ of
nonexpansive mappings on $H$ and $A$ is a strongly monotone and
lipschitzian mapping on $H$. 
Then we study convergence of the iterative sequence $\{x_n\}$ defined by 
$x_1 \in H$ and 
\begin{equation}\label{eqn:intro-hsdm}
 x_{n+1} = (I -  \lambda_n A ) T_n x_n 
\end{equation}
for $n\in \N$ in order to approximate the solution,
where $I$ is the identity mapping on $H$. 
If $A= I-u$ for some $u\in H$, 
then it is clear that $A$ is strongly monotone and lipschitzian, 
and~\eqref{eqn:intro-hsdm} is reduced to 
\begin{equation}\label{eqn:intro-halpern}
 x_{n+1} = \lambda_n u + (1- \lambda_n) T_n x_n. 
\end{equation}
We deal with these two types of iterations, and especially we focus on
the relationship between them; see \S\ref{sec:halpern-hsdm}. 

The iterative method defined by~\eqref{eqn:intro-hsdm} is called
the hybrid steepest descent method, which was introduced by
Yamada~\cite{MR1853237}.
He considered the variational inequality problem over the set
of common fixed points of a finite family of nonexpansive mappings and
proved strong convergence of the sequence generated by the method. 
We know many results by using the hybrid steepest descent method; 
see~\cites{%
NACA2009,
MR2504478,
MR2501545,
MR2028445,
MR2311630,
MR2418843,
MR2467157,
MR2588936,
MR2210870,
MR2418842,
MR2106274,
MR2536778,
MR2653741,
MR1977756,
MR2109044,
MR2338663,
MR2200594,
MR2303800,
MR2495919}.

The iterative method defined by~\eqref{eqn:intro-halpern} is called
the Halpern type iterative method; see Halpern~\cite{MR0218938}, 
Wittmann~\cite{MR1156581}, and Shioji and Takahashi~\cite{MR1415370}; 
see also~\cites{ICFPTA2009,AK,MR2338104}. 

\section{Preliminaries}\label{s:preliminaries}
Throughout the present paper, $H$ denotes a real Hilbert space
with the inner product $\ip{\,\cdot\,}{\,\cdot\,}$ and 
the norm $\norm{\,\cdot\,}$, 
$I$ the identity mapping on $H$, 
and $\N$ the set of positive integers. 

Let $C$ be a nonempty closed convex subset of $H$. 
A mapping $S\colon C\to C$ is said to be lipschitzian
if there exists a constant $\eta > 0$ such that 
$\norm{Sx - Sy} \leq \eta \norm{x-y}$ for all $x,y\in C$. 
In this case, $S$ is called an $\eta$-lipschitzian mapping. 
In particular, 
an $\eta$-lipschitzian mapping is said to be nonexpansive if $\eta=1$; 
an $\eta$-lipschitzian mapping is said to be a contraction if 
$0\leq \eta <1$. 
It is known that $\Fix(S)$ is closed and convex if $S$ is nonexpansive, 
where $\Fix(S)$ denotes the set of fixed points of $S$. 
The metric projection of $H$ onto $C$ is denoted by $P_C$
and we know that $P_C$ is nonexpansive. 
We also know the following; see~\cite{MR2548424}. 
\begin{lemma}
 \label{lm:projection}
 Let $x\in H$ and $z\in C$. 
 Then $z= P_C (x)$ if and only if $\ip{y- z}{x- z}\leq 0$
 for all $y\in C$. 
\end{lemma}

Let $C$ be a nonempty closed convex subset of $H$. 
Let $\{S_n\}$ be a sequence of self-mappings of $C$. 
We say that $\{S_n\}$ satisfies the condition~(Z) if every weak cluster
point of $\{x_n\}$ is a common fixed point of $\{S_n\}$ whenever 
$\{x_n \}$ is a bounded sequence in $C$ and $x_n - S_n x_n \to 0$; 
see~\cites{MR2671943,MR2377867,MR2529497,AKTT4,ICFPTA2009,ISBFS2009}.

A mapping $A\colon H\to H$ is said to be strongly monotone 
if there is a constant $\kappa > 0$ such that
$\ip{x-y}{Ax-Ay}\geq \kappa \norm{x-y}^2$ for all $x,y\in H$. 
In this case, $A$ is called a $\kappa$-strongly monotone mapping. 
In \S\ref{sec:halpern-hsdm}, 
we deal with the following variational inequality problem: 

\begin{problem}\label{p:VI}
 Let $\kappa$ and $\eta$ be positive real numbers such that 
 $\eta^2 < 2\kappa$. 
 Let $F$ a nonempty closed convex subset of $H$
 and $A\colon H\to H$ a $\kappa$-strongly monotone and
 $\eta$-lipschitzian mapping. 
 Then find $z\in F$ such that 
 \[
  \ip{y-z}{Az} \geq 0 \text{ for all } y \in F.
 \]
\end{problem}

The set of solution of Problem~\ref{p:VI} is denoted by $\VI(F,A)$. 
It is known that the solution set is a singleton; see
Lemma~\ref{lm:basic} below. 

\begin{remark}
 The assumption that $\eta^2 < 2\kappa$ in Problem~\ref{p:VI} is not
 restrictive. 
 Indeed, suppose that a $\kappa$-strongly monotone and
 $\eta$-lipschitzian mapping $A$ is given. 
 Let us choose a positive constant $\mu$ such that 
 $\mu < 2\kappa/ \eta^2$, 
 and define $\kappa' = \mu \kappa$ and $\eta'=\mu\eta$. 
 Then it is easy to verify that $(\eta')^2 < 2 \kappa'$, 
 $\mu A$ is $\kappa'$-strongly monotone and $\eta'$-lipschitzian, 
 and moreover, $\VI(F, A) = \VI(F, \mu A)$
 for every nonempty closed convex subset $F$ of $H$. 
\end{remark}

\begin{lemma}\label{lm:basic}
 Under the assumptions of Problem~\ref{p:VI}, the following hold:
\begin{enumerate}
 \item $\kappa\leq \eta$, $0 \leq 1- 2\kappa + \eta^2 <1$ and $I - A$
       is a $\theta$-contraction,
       where $\theta = \sqrt{1- 2\kappa + \eta^2}$. 
 \item Problem~\ref{p:VI} has a unique solution 
       and $\VI(F,A)= \Fix\bigl(P_F (I-A)\bigr)$. 
\end{enumerate}
\end{lemma}

\begin{proof}
 We first prove (1). 
 Since $A$ is $\kappa$-strongly monotone and $\eta$-lipschitzian, 
 it follows that 
 \[
 \kappa \norm{x-y}^2 \leq \ip{x-y}{Ax-Ay} 
 \leq \norm{x-y}\norm{Ax-Ay}\leq \eta \norm{x-y}^2 
 \]
 and 
 \begin{align*}
  \norm{(I- A)x - (I- A)y}^2
  &= \norm{x-y}^2 - 2 \ip{x-y}{Ax-Ay} + \norm{Ax- Ay}^2\\
  &\leq \norm{x-y}^2 - 2 \kappa \norm{x-y}^2 
  + \eta^2\norm{x- y}^2\\
  &= \left( 1 - 2\kappa + \eta^2 \right) \norm{x- y}^2
 \end{align*}
 for all $x,y\in H$. Therefore, $\kappa \leq \eta$. 
 By assumption, it is easy to check that $0\leq 1-2\kappa + \eta^2 <1$,
 and thus $I-A$ is a $\theta$-contraction. 

 We next prove (2). 
 Since $I-A$ is a contraction by~(1) and the metric projection $P_F$ is
 nonexpansive, $P_F (I-A)$ is a contraction on $H$. 
 The Banach contraction principle guarantees that $P_F(I-A)$ has a
 unique fixed point. 
 On the other hand, Lemma~\ref{lm:projection} shows that 
 $\VI(F,A)= \Fix\bigl(P_F(I-A) \bigr)$ and thus Problem~\ref{p:VI} has a
 unique solution. 
\end{proof}

We know the following result; see~\cites{MR2338104,NACA2009} and 
see also~\cites{MR2529497,ISBFS2009}. 

\begin{theorem}\label{thm:aktt}
 Let $H$ be a Hilbert space, $C$ a nonempty closed convex subset of $H$,
 $\{T_n\}$ a sequence of nonexpansive self-mappings of $C$
 with a common fixed point, $F$ the set of common fixed points of $\{T_n\}$, 
 and $\{\lambda_n\}$ a sequence in $[0,1]$ such that 
 \begin{equation}\label{eqn:lambda1}
 \lambda_n \to 0, \, 
 \sum_{n=1}^\infty \lambda_n =\infty
 \text{, and } 
 \sum_{n=1}^\infty \abs{\lambda_{n+1} - \lambda_n}<\infty.   
 \end{equation}
 Suppose that $\{T_n\}$ satisfies the condition~(Z) and
 \begin{equation}\label{eqn:aktt}
  \sum_{n=1}^\infty \sup\{\norm{T_{n+1} y - T_n y}: y\in D\} < \infty
 \end{equation}
 for every nonempty bounded subset $D$ of $C$. 
 Let $x,u$ be points in $C$ and $\{x_n\}$ a sequence 
 defined by $x_1 = x$ and \eqref{eqn:intro-halpern}
 for $n\in \N$. 
 Then $\{ x_n \}$ converges strongly to $P_F (u)$. 
\end{theorem}

We also know the following result; see~\cites{ICFPTA2009,AK}. 

\begin{theorem}\label{thm:ak}
 Let $H$ be a Hilbert space, $C$ a nonempty closed convex subset of $H$,
 $\{S_n\}$ a sequence of nonexpansive self-mappings of $C$
 with a common fixed point, $F$ the set of common fixed points of
 $\{S_n\}$. 
 Let $\{\lambda_n\}$ and $\{\beta_n\}$ be a sequences in $[0,1]$ such
 that 
 \begin{equation*}
 \lambda_n \to 0,\,
 \sum_{n=1}^\infty \lambda_n =\infty\text{, and }
 0 < \liminf_{n\to\infty} \beta_n \leq \limsup_{n\to\infty} \beta_n <1. 
 \end{equation*}
 Suppose that $\{S_n\}$ satisfies the condition~(Z) and
 \begin{equation}\label{eqn:R}
  \sup\{\norm{S_{n+1} y - S_n y}: y\in D\} \to 0
 \end{equation}
 for every nonempty bounded subset $D$ of $C$. 
 Let $x,u$ be points in $C$ and $\{x_n\}$ a sequence defined by 
 $x_1 = x$ and 
 \[
 x_{n+1} = \lambda_n u + (1- \lambda_n) 
 \bigl( \beta_n x_n + (1 - \beta_n)S_n x_n \bigr) 
 \]
 for $n\in \N$. 
 Then $\{ x_n \}$ converges strongly to $P_F (u)$. 
\end{theorem}

The following lemma is well known; see~\cites{%
MR1353071,MR1086345,MR1889384,MR1911872,MR2338104}. 

\begin{lemma}\label{lm:seq}
 Let $\{ \epsilon_n \}$ be a sequence of nonnegative real numbers, 
 $\{ \gamma_n \}$ a sequence of real numbers, 
 and $\{ \lambda_n \}$ a sequence in $[0,1]$. 
 Suppose that 
 $\epsilon_{n+1} \leq (1- \lambda_n) \epsilon_n + \lambda_n \gamma_n$
 for every $n\in\N$, $\limsup_{n\to\infty} \gamma_n \leq 0 $, 
 and $\sum_{n=1}^\infty \lambda_n = \infty$. Then $\epsilon_n \to 0$. 
\end{lemma}

\section{Convergence theorems by the hybrid steepest descent method}
\label{sec:halpern-hsdm}
In this section, we deal with the variational inequality problem 
over the set of common fixed points of a sequence of nonexpansive
mappings; see Problem~\ref{p:VI-FT-A} below. 
We first investigate the relationship between the hybrid steepest
descent method and the Halpern type iterative method
(Theorem~\ref{thm:halpern2hsdm}). 
Then, by using Theorems~\ref{thm:aktt} and~\ref{thm:ak}, 
we show two convergence theorems 
by the hybrid steepest descent method for this problem. 

\begin{problem}\label{p:VI-FT-A}
 Let $H$ be a Hilbert space, $\{T_n\}$ a sequence of nonexpansive 
 self-mappings of $H$ with a common fixed point, $F$ the set of common
 fixed points of $\{T_n\}$, and $A\colon H\to H$ a $\kappa$-strongly
 monotone and $\eta$-lipschitzian mapping,
 where $\kappa$ and $\eta$ are positive real numbers such that 
 $\eta^2 < 2\kappa$. 
 Then find $z\in F$ such that 
 \[
  \ip{y-z}{Az} \geq 0 \text{ for all } y \in F.
 \]
\end{problem}

Using the technique in \cite{MR2273529}, 
we can prove the following theorem, 
which shows that every sequence generated by the hybrid steepest descent
method for Problem~\ref{p:VI-FT-A} is convergent whenever so is every
sequence generated by the Halpern type iterative method for the sequence
of nonexpansive mappings. 

\begin{theorem}\label{thm:halpern2hsdm}
 Let $H$, $\{T_n\}$, $F$, $\kappa$, $\eta$, and $A$ be the same as in
 Problem~\ref{p:VI-FT-A}. 
 Let $\{\lambda_n\}$ be a sequence in $[0,1]$ such that 
 $\sum_{n=1}^\infty \lambda_n = \infty$. 
 Suppose that for any $(x,u)\in H\times H$, 
 the sequence $\{x_n\}$ defined by $x_1 = x$ and 
 \begin{equation}\label{eqn:I-u}
  x_{n+1} = \lambda_n u + (1-\lambda_n) T_n x_n  
 \end{equation}
  for $n\in \N$ converges strongly to $P_F(u)$. 
 Let $y$ be a point in $H$ and $\{y_n\}$ a sequence defined by
 $y_1 = y$ and
 \begin{equation}\label{eqn:hsdm}
  y_{n+1} = (I - \lambda_n A)T_n y_n  
 \end{equation}
 for $n\in \N$. 
 Then $\{y_n\}$ converges strongly to the unique solution of
 Problem~\ref{p:VI-FT-A}.
\end{theorem}

\begin{proof}
 Set $f_n = (I-A)T_n$ for $n\in \N$. 
 Since $T_n$ is nonexpansive, $f_n$ is a $\theta$-contraction on $H$ by
 Lemma~\ref{lm:basic}, where $\theta=\sqrt{1-2\kappa + \eta^2}$.
 Let $w$ be the fixed point of $P_F \circ f_1$ and 
 $\{x_n \}$ a sequence defined by $x_1 = y$ and 
 \[
  x_{n+1} = \lambda_n f_1 (w) + (1 - \lambda_n) T_n x_n   
 \]
 for $n\in\N$. Then 
 \begin{equation} \label{eqn:y_n2v}
   x_n \to P_F \bigl( f_1 (w) \bigr) = w
 \end{equation}
 by assumption. Since $T_n$ is nonexpansive and $f_n$ is a
 $\theta$-contraction, it follows from $f_1(w) = f_n(w)$ that
 \begin{align*}
  \norm{x_{n+1} - y_{n+1}}
  &= \norm{(1-\lambda_n) (T_n x_n - T_n y_n)
  + \lambda_n \bigl( f_1 (w) - f_n (y_n) \bigr)}\\
  &\leq (1-\lambda_n)\norm{T_n x_n - T_n y_n}
  + \lambda_n \norm{f_n (w) - f_n (y_n)}\\
  &\leq (1-\lambda_n)\norm{x_n - y_n} + \lambda_n \theta \norm{w- y_n}\\
  &\leq (1-\lambda_n)\norm{x_n - y_n} + \lambda_n \theta (
  \norm{w - x_n} +  \norm{x_n - y_n}) \\
  &\leq \bigl( 1 - (1-\theta) \lambda_n \bigr) \norm{x_n - y_n} + 
  (1-\theta)\lambda_n \frac{\theta}{1- \theta} \norm{x_n - w}
 \end{align*}
 for every $n\in\N$. 
 Since $\sum_{n=1}^\infty (1- \theta)\lambda_n = \infty$, 
 \eqref{eqn:y_n2v} and Lemma~\ref{lm:seq} show that 
 $x_n - y_n \to 0$. Therefore, we conclude that $\{y_n\}$ converges
 strongly to $w = P_F \bigl( (I-A)T_1 w\bigr) = P_F (I-A)w$, 
 which is the unique solution of Problem~\ref{p:VI-FT-A}
 by Lemma~\ref{lm:basic}. 
 This completes the proof. 
\end{proof}

Using Theorems~\ref{thm:aktt} and~\ref{thm:halpern2hsdm}, 
we obtain the following:

\begin{theorem}
 Let $H$, $\{T_n\}$, $F$, $\kappa$, $\eta$, and $A$ be the same as in
 Problem~\ref{p:VI-FT-A}. 
 Let $\{\lambda_n\}$ be a sequence in $[0,1]$ such that
 \eqref{eqn:lambda1} holds. 
 Suppose that $\{T_n\}$ satisfies the condition~(Z) and
 \eqref{eqn:aktt} holds 
 for every nonempty bounded subset $D$ of $C$. 
 Let $y$ be a point in $H$ and $\{y_n\}$ a sequence defined by 
 $y_1 = y$ and~\eqref{eqn:hsdm} for $n\in \N$. 
 Then $\{y_n\}$ converges strongly to the unique solution of
 Problem~\ref{p:VI-FT-A}.  
\end{theorem}

\begin{proof}
 Let $(x,u)\in H\times H$ be fixed. 
 Then it follows from Theorem~\ref{thm:aktt} that 
 the sequence $\{x_n\}$ defined by $x_1 = x$ and \eqref{eqn:I-u} for
 $n\in \N$ converges strongly to $P_F(u)$. 
 Therefore, Theorem~\ref{thm:halpern2hsdm} implies the conclusion. 
\end{proof}

Using Theorem~\ref{thm:halpern2hsdm} and other known results, 
we also obtain the following: 

\begin{theorem}[Iemoto and Takahashi~\cite{MR2501545}*{Theorem 3.1}]
 Let $H$, $\{T_n\}$, $F$, $\kappa$, $\eta$, and $A$ be the same as in
 Problem~\ref{p:VI-FT-A}. 
 Let $\{\lambda_n\}$ be a sequence in $[0,1]$ such that 
 \[
  \lambda_n \to 0 \text{ and } \sum_{n=1}^\infty \lambda_n = \infty
 \]
 and $\{\gamma_n\}$ a sequence in $[a,b]$, where $0 < a \leq b < 1$. 
 For each $n\in \N$ and $k\in \{ 1,2, \dotsc, n+1 \}$, 
 define a mapping $U_{n,k}$ by  $U_{n,n+1} = I$ and 
 \[
  U_{n,k}= \gamma_k T_k U_{n,k+1} + (1-\gamma_k) I. 
 \]
 Let $y$ be a point in $H$ and $\{ y_n \}$ a sequence defined by 
 $y_1 = y$ and
 \begin{equation}\label{eqn:iemoto-T}
  y_{n+1} = ( I- \lambda_n A ) U_{n,1} y_n 
 \end{equation}
 for $n\in \N$. 
 Then $\{y_n\}$ converges strongly to the unique solution 
 of Problem~\ref{p:VI-FT-A}. 
\end{theorem}

\begin{proof}
 Set $S_n = T_1 U_{n,2}$ for $n \in \N$.
 Then it is clear that each $S_n$ is nonexpansive. 
 It is known that 
 \[
 \Fix(S_n) = \Fix(U_{n,1})  = \bigcap_{k=1}^n \Fix(T_k)
 \]
 by \cite{MR1800669}*{Lemma 3.2}; see also \cite{AKTT4}*{Lemma 4.2}. 
 Hence we have 
 \[
 \bigcap_{n=1}^\infty \Fix(S_n) 
 = \bigcap_{n=1}^\infty \Fix(U_{n,1})
 = \bigcap_{n=1}^\infty \bigcap_{k=1}^n \Fix(T_k) = F. 
 \]
 It is also known that $\{S_n\}$ satisfies the condition~(Z) 
 and \eqref{eqn:R} holds for every nonempty bounded subset $D$ of $H$;
 see~\cite{AKTT4}, \cite{MR2529497}, \cite{ISBFS2009}, and~\cite{AK}. 
 Thus, for any $(x,u)\in H\times H$, 
 it follows from Theorem~\ref{thm:ak} that 
 the sequence $\{x_n\}$ defined by $x_1 = x$ and 
 \[
 x_{n+1} = \lambda_n u + (1-\lambda_n) U_{n,1} x_n 
 = \lambda_n u + (1-\lambda_n) 
 \bigl( (1-\gamma_1) x_n + \gamma_1 S_n x_n \bigr)
 \]
 for $n\in \N$ converges strongly to $P_F(u)$. 
 Therefore, Theorem~\ref{thm:halpern2hsdm} implies the conclusion. 
\end{proof}


\begin{bibdiv}
 \begin{biblist}
\bib{ICFPTA2009}{article}{
   author={Aoyama, Koji},
   title={An iterative method for fixed point problems for sequences of
  nonexpansive mappings}, 
   conference={
      title={Fixed Point theory and its Applications},
   },
   book={
      publisher={Yokohama Publ., Yokohama},
   },
  date={2010},
  pages={1--7},
}

\bib{NACA2009}{article}{
   author={Aoyama, Koji},
   title={An iterative method for a variational inequality problem over the
  common fixed point set of nonexpansive mappings},
   conference={
      title={Nonlinear analysis and convex analysis},
   },
   book={
      publisher={Yokohama Publ., Yokohama},
   },
  date={2010},
  pages={21--28},
}

\bib{ISBFS2009}{article}{
   author={Aoyama, Koji},
   title={Asymptotic fixed points of sequences of quasi-nonexpansive
  type mappings}, 
  journal={Proceedings of the 3rd International Symposium on Banach
  and Function Spaces},
  status={to appear},
}

\bib{AK}{article}{
   author={Aoyama, Koji},
   author={Kimura, Yasunori},
   title={Strong convergence theorems for strongly nonexpansive
  sequences},
  status={submitted},
}

\bib{MR2338104}{article}{
   author={Aoyama, Koji},
   author={Kimura, Yasunori},
   author={Takahashi, Wataru},
   author={Toyoda, Masashi},
   title={Approximation of common fixed points of a countable family of
   nonexpansive mappings in a Banach space},
   journal={Nonlinear Anal.},
   volume={67},
   date={2007},
   pages={2350--2360},
}

\bib{MR2377867}{article}{
   author={Aoyama, Koji},
   author={Kimura, Yasunori},
   author={Takahashi, Wataru},
   author={Toyoda, Masashi},
   title={On a strongly nonexpansive sequence in Hilbert spaces},
   journal={J. Nonlinear Convex Anal.},
   volume={8},
   date={2007},
   pages={471--489},
}
		
\bib{AKTT4}{article}{
   author={Aoyama, Koji},
   author={Kimura, Yasunori},
   author={Takahashi, Wataru},
   author={Toyoda, Masashi},
   title={Strongly nonexpansive sequences and their applications in
 Banach spaces},
   conference={
      title={Fixed Point theory and its Applications},
   },
   book={
      publisher={Yokohama Publ., Yokohama},
   },
   date={2008},
   pages={1--18},
}

\bib{MR2671943}{article}{
   author={Aoyama, Koji},
   author={Kohsaka, Fumiaki},
   author={Takahashi, Wataru},
   title={Shrinking projection methods for firmly nonexpansive mappings},
   journal={Nonlinear Anal.},
   volume={71},
   date={2009},
   pages={e1626--e1632},

}

\bib{MR2529497}{article}{
   author={Aoyama, Koji},
   author={Kohsaka, Fumiaki},
   author={Takahashi, Wataru},
   title={Strongly relatively nonexpansive sequences in Banach spaces and
   applications},
   journal={J. Fixed Point Theory Appl.},
   volume={5},
   date={2009},
   pages={201--224},
}

\bib{MR2418843}{article}{
   author={Ceng, Lu-Chuan},
   author={Xu, Hong-Kun},
   author={Yao, Jen-Chih},
   title={A hybrid steepest-descent method for variational inequalities in
   Hilbert spaces},
   journal={Appl. Anal.},
   volume={87},
   date={2008},
   pages={575--589},
}

\bib{MR2467157}{article}{
   author={Ceng, Lu-Chuan},
   author={Ansari, Qamrul Hasan},
   author={Yao, Jen-Chih},
   title={Mann-type steepest-descent and modified hybrid steepest-descent
   methods for variational inequalities in Banach spaces},
   journal={Numer. Funct. Anal. Optim.},
   volume={29},
   date={2008},
   pages={987--1033},
}

\bib{MR2504478}{book}{
   author={Chidume, Charles},
   title={Geometric properties of Banach spaces and nonlinear iterations},
   series={Lecture Notes in Mathematics},
   volume={1965},
   publisher={Springer-Verlag London Ltd.},
   place={London},
   date={2009},
   pages={xviii+326},
}

\bib{MR0218938}{article}{
 author={Halpern, Benjamin},
 title={Fixed points of nonexpanding maps},
 journal={Bull. Amer. Math. Soc.},
 volume={73},
 date={1967},
 pages={957--961},
}

\bib{MR2501545}{article}{
   author={Iemoto, Shigeru},
   author={Takahashi, Wataru},
   title={Strong convergence theorems by a hybrid steepest descent method
   for countable nonexpansive mappings in Hilbert spaces},
   journal={Sci. Math. Jpn.},
   volume={69},
   date={2009},
   pages={227--240},
}

\bib{MR2588936}{article}{
   author={Kirihara, Ikuo},
   author={Kurokawa, Yu},
   author={Takahashi, Wataru},
   title={Strong convergence theorem for quadratic minimization problem with
   countable constraints},
   journal={J. Nonlinear Convex Anal.},
   volume={10},
   date={2009},
   pages={383--393},
}

\bib{MR1353071}{article}{
   author={Liu, Li Shan},
   title={Ishikawa and Mann iterative process with errors for nonlinear
   strongly accretive mappings in Banach spaces},
   journal={J. Math. Anal. Appl.},
   volume={194},
   date={1995},
   pages={114--125},
}

\bib{MR2210870}{article}{
   author={Marino, Giuseppe},
   author={Xu, Hong-Kun},
   title={A general iterative method for nonexpansive mappings in Hilbert
   spaces},
   journal={J. Math. Anal. Appl.},
   volume={318},
   date={2006},
   pages={43--52},
}

\bib{MR2418842}{article}{
   author={Miao, Yu},
   author={Song, Yisheng},
   title={Weak and strong convergence of a new scheme for two non-expansive
   mappings in Hilbert spaces},
   journal={Appl. Anal.},
   volume={87},
   date={2008},
   pages={567--574},
}


\bib{MR2106274}{article}{
   author={O'Hara, John G.},
   author={Pillay, Paranjothi},
   author={Xu, Hong-Kun},
   title={Iterative approaches to convex minimization problems},
   journal={Numer. Funct. Anal. Optim.},
   volume={25},
   date={2004},
   pages={531--546},
}

\bib{MR2536778}{article}{
   author={Qin, Xiaolong},
   author={Cho, Yeol Je},
   author={Kang, Shin Min},
   title={Some results on non-expansive mappings and relaxed cocoercive
   mappings in Hilbert spaces},
   journal={Appl. Anal.},
   volume={88},
   date={2009},
   pages={1--13},
}

\bib{MR1415370}{article}{
   author={Shioji, Naoki},
   author={Takahashi, Wataru},
   title={Strong convergence of approximated sequences for nonexpansive
   mappings in Banach spaces},
   journal={Proc. Amer. Math. Soc.},
   volume={125},
   date={1997},
   pages={3641--3645},
}

\bib{MR2273529}{article}{
   author={Suzuki, Tomonari},
   title={Moudafi's viscosity approximations with Meir-Keeler contractions},
   journal={J. Math. Anal. Appl.},
   volume={325},
   date={2007},
   pages={342--352},
}

\bib{MR2548424}{book}{
   author={Takahashi, Wataru},
   title={Introduction to nonlinear and convex analysis},
   publisher={Yokohama Publ., Yokohama},
   date={2009},
   pages={iv+234},
}

\bib{MR1800669}{article}{
   author={Takahashi, W.},
   author={Shimoji, K.},
   title={Convergence theorems for nonexpansive mappings and feasibility
   problems},
   journal={Math. Comput. Modelling},
   volume={32},
   date={2000},
   pages={1463--1471},
}

\bib{MR2653741}{article}{
   author={Tian, Ming},
   title={A general iterative algorithm for nonexpansive mappings in Hilbert
   spaces},
   journal={Nonlinear Anal.},
   volume={73},
   date={2010},
   pages={689--694},
}

\bib{MR1086345}{article}{
   author={Weng, Xinlong},
   title={Fixed point iteration for local strictly pseudo-contractive
   mapping},
   journal={Proc. Amer. Math. Soc.},
   volume={113},
   date={1991},
   pages={727--731},
} 

\bib{MR1156581}{article}{
   author={Wittmann, Rainer},
   title={Approximation of fixed points of nonexpansive mappings},
   journal={Arch. Math. (Basel)},
   volume={58},
   date={1992},
   pages={486--491},
}

\bib{MR1889384}{article}{
   author={Xu, Hong-Kun},
   title={Another control condition in an iterative method for nonexpansive
   mappings},
   journal={Bull. Austral. Math. Soc.},
   volume={65},
   date={2002},
   pages={109--113},

}

\bib{MR1911872}{article}{
   author={Xu, Hong-Kun},
   title={Iterative algorithms for nonlinear operators},
   journal={J. London Math. Soc. (2)},
   volume={66},
   date={2002},
   pages={240--256},
}

\bib{MR1977756}{article}{
   author={Xu, H. K.},
   title={An iterative approach to quadratic optimization},
   journal={J. Optim. Theory Appl.},
   volume={116},
   date={2003},
   pages={659--678},
}


\bib{MR2028445}{article}{
   author={Xu, H. K.},
   author={Kim, T. H.},
   title={Convergence of hybrid steepest-descent methods for variational
   inequalities},
   journal={J. Optim. Theory Appl.},
   volume={119},
   date={2003},
   pages={185--201},
}

\bib{MR1853237}{article}{
   author={Yamada, Isao},
   title={The hybrid steepest descent method for the variational inequality
   problem over the intersection of fixed point sets of nonexpansive
   mappings},
   conference={
      title={Inherently parallel algorithms in feasibility and
  optimization their applications},
      address={Haifa},
      date={2000},
   },
   book={
      series={Stud. Comput. Math.},
      volume={8},
      publisher={North-Holland},
      place={Amsterdam},
   },
   date={2001},
   pages={473--504},
}

\bib{MR2109044}{article}{
   author={Yamada, Isao},
   author={Ogura, Nobuhiko},
   title={Hybrid steepest descent method for variational inequality problem
   over the fixed point set of certain quasi-nonexpansive mappings},
   journal={Numer. Funct. Anal. Optim.},
   volume={25},
   date={2004},
   pages={619--655},
}


\bib{MR2311630}{article}{
   author={Yao, Yonghong},
   title={A general iterative method for a finite family of nonexpansive
   mappings},
   journal={Nonlinear Anal.},
   volume={66},
   date={2007},
   pages={2676--2687},
}

\bib{MR2338663}{article}{
   author={Yao, Yonghong},
   author={Noor, Muhammad Aslam},
   title={On modified hybrid steepest-descent methods for general
   variational inequalities},
   journal={J. Math. Anal. Appl.},
   volume={334},
   date={2007},
   pages={1276--1289},
}

\bib{MR2200594}{article}{
   author={Zeng, Liu Chuan},
   author={Wong, N. C.},
   author={Yao, J. C.},
   title={Convergence of hybrid steepest-descent methods for generalized
   variational inequalities},
   journal={Acta Math. Sin. (Engl. Ser.)},
   volume={22},
   date={2006},
   pages={1--12},
}

\bib{MR2303800}{article}{
   author={Zeng, L. C.},
   author={Wong, N. C.},
   author={Yao, J. C.},
   title={Convergence analysis of modified hybrid steepest-descent methods
   with variable parameters for variational inequalities},
   journal={J. Optim. Theory Appl.},
   volume={132},
   date={2007},
   pages={51--69},
}

\bib{MR2495919}{article}{
   author={Zeng, L. C.},
   author={Schaible, S.},
   author={Yao, J. C.},
   title={Hybrid steepest descent methods for zeros of nonlinear operators
   with applications to variational inequalities},
   journal={J. Optim. Theory Appl.},
   volume={141},
   date={2009},
   pages={75--91},
}
 \end{biblist}
\end{bibdiv}
\end{document}